\documentclass[prb,preprint,eqsecnum]{revtex4}%
\usepackage{amsmath}
\usepackage{graphicx}
\usepackage{amssymb}
\usepackage{graphicx}
\usepackage{amsfonts}%
\providecommand{\U}[1]{\protect\rule{.1in}{.1in}}

\begin{document}
\title{A new integral representation for the Riemann Zeta function}
\author{Sandeep Tyagi, Christian Holm}
\affiliation{Frankfurt Institute for Advanced Studies, Max-von-Laue-Stra{\ss }e 1,}
\affiliation{D-60438 Frankfurt am Main, Germany}

\begin{abstract}
A new integral representation for the Riemann zeta function is derived. This
new representation of $\zeta\left( s\right) $ covers the important region of
the complex plane where $0<\operatorname{Re}\left[  s\right]  <1$. Using this
new representation, we obtain new functional identities for the Riemann zeta function.

\end{abstract}
\maketitle


\section{Introduction}

The Riemann zeta function has an illustrious history of 150 years
\cite{edwards}. It arises in many branches of mathematics and, and is also
intimately related to the distribution of prime numbers. The most facinating
thing related to the Riemann zeta function is the Riemann's hypothesis which
states that all non-trivial zeros of the Riemann zeta function $\zeta\left(
s\right)  $ fall on the strip $\operatorname{Re}[s]=1/2.$ This hypothesis has
so far not been proven despite the efforts of some of the most brilliant
scientists in mathematics such as Riemann, Hilbert, Hardy and Polya
\cite{edwards}. There are a number of integral representations of the Riemann
zeta function. A few of them cover the important region of $0\leq$
Re$[s]\leq1$. In this paper, we report a new integral represntation.

The Riemann zeta function is defined as
\begin{equation}
\zeta\left(  s\right)  =\sum_{n=1}^{\infty}\frac{1}{n^{s}}\nonumber
\end{equation}
However, this definition is valid only in the region $\operatorname{Re}[s]>1$.
By analytical continuation, the definition can be extended to the whole
complex plane. The new formula valid for the whole complex plane is%
\begin{equation}
\zeta\left(  s\right)  =-\frac{1}{1-2^{\left(  1-s\right)  }}\sum
_{n=1}^{\infty}\frac{\left(  -1\right)  ^{n}}{n^{s}}.\label{z2}%
\end{equation}
The Riemann zeta function has just a single simple pole in the whole complex
plane which lies at $s=1$. The $\zeta\left(  s\right)  $ has a number of
integral representations; a comprehensive list can be obtained from the
mathworld website \cite{wolfram}. Our aim in this paper is to derive a new
integral representation and obtain new functional identity using it.

\section{A new representation}

We obtain this representation by rewriting a general term of Eq.~(\ref{z2}) in
an integral form and then taking the sum over $n$ inside the integral sign. We
start with Eq.~(\ref{z2}). We will first work with a special case of
$s=1/2+$i$b$. In the end, we will use analytical continuation to obtain the
result for any $0<\operatorname{Re}[s]<2$. For $s=1/2+\text{i}b$ we have from
Eq.~(\ref{z2}),%
\begin{equation}
\zeta\left(  1/2+\text{i}b\right)  =-\frac{1}{1-2^{\left(  1/2-\text{i}%
b\right)  }}\sum_{n=1}^{\infty}\left(  -1\right)  ^{n}\frac{n^{-\text{i}b}%
}{\sqrt{n}}.\label{half}%
\end{equation}
At this point, we cannot proceeed any further as the summation over $n$ cannot
be performed in a straightforward way. We have to re-express the term
$n^{-\text{i}b}/\sqrt{n}$ in Eq.~(\ref{half}) so that the summation over $n$
can be carried out. To this end, we would like to have an integral involving
the dummy index $x$ and the parameters $n$ and $b$ such that this integral is
equal to $n^{-\text{i}b}/\sqrt{n}$ apart from some factors independent of $n$.
To achieve such an integral, we can try the complex residue theorem. It
dictates that the denominator should have a form of $x^{2}+n$ or $x^{4}+n^{2}$
because then the residues would give rise to terms having $\sqrt{n}$. One can
go ahead and try it with $x^{2}+n$, but that leads to difficulties in the
series summation later on. So, let us try with $x^{4}+n^{2}$, which does lead
to a sensible series. To get rid of some other unwanted terms that would arise
out of the presence of residues, we consider the integral%
\begin{equation}
\int_{-\infty}^{+\infty}\frac{x^{2}x^{-2\text{i}b}}{x^{4}+n^{2}}dx.\label{int}%
\end{equation}
This integral has been chosen such that when $x^{2}$ is substituted by
$\sqrt{n}$, then the $x^{2}$ in the numerator will be cancelled by a term
arising in the denominator and the $x^{-2\text{i}b}$ would give rise to
$n^{-\text{i}b}$. Finally a term of $\sqrt{n}$ would arise as a result of the
other residue. The integral in Eq.~(\ref{int1}) can be easily evaluated using
the complex residue theorem by choosing an appropriate contour. Even though we
have assumed $b$ to be real when writing the integral in Eq.~(\ref{int1}), it
turns out that the integral is defined for all $b$ in the range of
$-3/2<\operatorname{Im}\left[  b\right]  <1/2$. The integral can be performed
to yield
\begin{equation}
\int_{-\infty}^{+\infty}\frac{x^{2}x^{-2\text{i}b}}{x^{4}+n^{2}}dx=\frac{\pi
}{2}\frac{n^{-\text{i}b}}{\sqrt{n}}\csc\left(  \frac{\pi}{2}\left(
1/2+\text{i}b\right)  \right)  \quad\quad\quad-3/2<\operatorname{Im}\left[
b\right]  <1/2.\label{int1}%
\end{equation}
Combining Eq.~(\ref{half}) and (\ref{int1}) we obtain%
\begin{equation}
\zeta\left(  1/2+\text{i}b\right)  =-\frac{2}{\pi}\frac{\sin\left(  \frac{\pi
}{2}\left(  1/2+\text{i}b\right)  \right)  }{\left(  1-2^{\left(
1/2-\text{i}b\right)  }\right)  }\sum_{n=1}^{\infty}\left(  -1\right)
^{n}\int_{-\infty}^{+\infty}\frac{x^{2}x^{-2\text{i}b}}{x^{4}+n^{2}}dx.
\end{equation}
The transformation above has led to a form where the summation over $n$ can
now be performed using well known results. We can interchange the summation
and integration and use the summation formula,%
\begin{equation}
\sum_{n=1}^{\infty}\frac{\left(  -1\right)  ^{n}}{x^{4}+n^{2}}=\frac{\pi
x^{2}\csc\left(  \pi x^{2}\right)  -1}{2x^{4}},
\end{equation}
to write%
\begin{equation}
\zeta\left(  1/2+\text{i}b\right)  =\frac{2}{\left(  1-2^{\left(
1/2-\text{i}b\right)  }\right)  }\sin\left(  \frac{\pi}{2}\left(
1/2+\text{i}b\right)  \right)  \int_{0}^{+\infty}\frac{x^{-2\text{i}b}}{\pi
x^{2}}\left(  1-\pi x^{2}\csc\left(  \pi x^{2}\right)  \right)  dx.\label{f1}%
\end{equation}
Note that the integral is symmetrical and for this reason we have written%
\[
\int_{-\infty}^{+\infty}...dx=2\int_{0}^{+\infty}...dx.
\]
Eq.~(\ref{f1}) is our desired integral representation. It is defined only for
$\operatorname{Re}\left[  s\right]  =1/2$, although the integral in
Eq.~(\ref{int}) is valid for $-3/2<\operatorname{Im}\left[  b\right]  <1/2$.
Writing $b=b^{\prime}+\text{i}q$ where $-3/2<q<1/2$ and $b^{\prime}$ is real,
we can express%
\[
1/2+\text{i}b=(1/2-q)+\text{i}b^{\prime}.
\]
Thus the real part of $s=1/2+\text{i}b$ actually changes between $0$ and $2$
as $q$ changes varies from $-3/2$ and $1/2.$ Substituting $s=1/2+\text{i}b$ in
Eq.~(\ref{f1}) we obtain%
\begin{equation}
\zeta\left(  s\right)  =\frac{2}{\left(  1-2^{1-s}\right)  }\sin\left(
\frac{\pi s}{2}\right)  \int_{0}^{\infty}\frac{x^{-2s}}{\pi x}\left(  1-\pi
x^{2}\csc\left[  \pi x^{2}\right]  \right)  dx\quad\quad\quad
0<\operatorname{Re}[s]<2.\label{f2}%
\end{equation}
Changing variable $\pi x^{2}=y$ we can write%
\begin{equation}
\zeta\left(  s\right)  =\frac{\pi^{s-1}}{\left(  1-2^{1-s}\right)  }%
\sin\left(  \frac{\pi s}{2}\right)  \int_{0}^{\infty}y^{-s}\left(  \frac{1}%
{y}-\frac{1}{\sinh\left(  y\right)  }\right)  dy\text{ }\quad\quad
\quad0<\operatorname{Re}\left[  s\right]  <2.\label{f3}%
\end{equation}
Note that Eqs.~(\ref{f2}) and (\ref{f3}) are valid for $0<\operatorname{Re}%
\left[  s\right]  <2.$ However, the Riemann zeta function has a simple pole on
$s=1$. Thus, these equations make sense only when $\operatorname{Im}[s]\neq0$
for $\operatorname{Re}[s]$ $=1$.

There are two formulas in the literature that resemble our formula in
Eq.~(\ref{f3}). The first one is
\begin{align}
\zeta\left(  s\right)   &  =\frac{2^{s-1}}{\left(  1-2^{1-s}\right)
\Gamma\left(  s\right)  }\int_{0}^{\infty}y^{s-1}\frac{\exp\left(  -y\right)
}{\cosh\left(  y\right)  }dy & \quad\quad\quad\operatorname{Re}\left[
s\right]   &  >0~\label{f4}\\
&  =\frac{2^{s-1}}{\left(  1-2^{1-s}\right)  \Gamma\left(  s\right)  }\int
_{0}^{\infty}\frac{y^{s-1}}{\exp\left(  y\right)  +1}dy & \quad\quad
\quad\operatorname{Re}\left[  s\right]   &  >0, \label{f5}%
\end{align}
where $\Gamma$ stands for the Gamma function. The second formula is due to
Ramanujan \cite{edwards} and it is given by%
\begin{equation}
\zeta\left(  s\right)  =\frac{\sin\left(  \pi s\right)  }{\pi}\int_{0}%
^{\infty}y^{-s}\left(  \log\left(  y\right)  -\frac{\Gamma^{\prime}\left(
y\right)  }{\Gamma\left(  y\right)  }\right)  dy\quad\quad\quad0\leq
\operatorname{Re}\left[  s\right]  <1. \label{f6}%
\end{equation}

There is another way to view \ref{f3}. The $\sinh(y)$ in denominstor will lead
to poles at integers with residues $(-1)^{n}.$ The presence of $1/y$ is
basically to remove the pole at $n=0.$ Thus \ref{f3} can be derived using
residue theorem as well.

\section{New functional identity}

So far there have been only a few functional forms and identities satisfied by
the Riemann zeta. One of the most important functional relation is given by%
\begin{equation}
\Gamma\left(  \frac{s}{2}-1\right)  \pi^{-\frac{s}{2}}\zeta\left(  s\right)
=\Gamma\left(  \frac{1-s}{2}-1\right)  \pi^{-\frac{1-s}{2}}\zeta\left(
1-s\right)  .
\end{equation}
Our representation leads to other functional identities involving the Riemann
zeta and Gamma functions. To derive a new functional representation we use the
integral formula%
\begin{equation}
\int_{0}^{\infty}\frac{y^{k}}{\sinh\left(  y\right)  }dy=\left(
2-2^{-k}\right)  \Gamma\left(  1+k\right)  \zeta\left(  1+k\right)  \quad
\quad\quad\text{where }k>0. \label{intf}%
\end{equation}
We can transform our general representation given in Eq.~(\ref{f3}) into this
form by using the series expansion
\begin{equation}
\sinh\left(  y\right)  -y=\frac{y^{3}}{3!}+\frac{y^{5}}{5!}+\cdots.
\end{equation}
This allows us to write Eq.~(\ref{f3}) as
\begin{align}
\zeta\left(  s\right)   &  =\frac{\pi^{s-1}}{\left(  1-2^{1-s}\right)  }%
\sin\left(  \frac{\pi s}{2}\right)  \int_{0}^{\infty}y^{-s-1}\left(
\frac{\sinh\left(  y\right)  -y}{\sinh\left(  y\right)  }\right)
dy\nonumber\\
&  =\text{ }\frac{\pi^{s-1}}{\left(  1-2^{1-s}\right)  }\sin\left(  \frac{\pi
s}{2}\right)  \sum_{n=1}^{\infty}\frac{1}{\left(  2n+1\right)  !}\int
_{0}^{\infty}y^{-s-1}\left(  \frac{y^{\left(  2n+1\right)  }}{\sinh\left(
y\right)  }\right)  dy.
\end{align}
Inserting the integral formula of Eq.~(\ref{intf}) leads to the following new
functional representation of the Riemann zeta function:%
\begin{equation}
\zeta\left(  s\right)  =\frac{\pi^{s-1}}{\left(  1-2^{1-s}\right)  }%
\sin\left(  \frac{\pi s}{2}\right)  \sum_{n=1}^{\infty}\frac{\left(
2-2^{-2n+s}\right)  }{\left(  2n+1\right)  !}\Gamma\left(  2n-s+1\right)
\zeta\left(  2n-s+1\right)  . \label{newf}%
\end{equation}
A part of the series in Eq.~(\ref{newf}) can be summed up analytically:
\begin{align}
\sum_{n=1}^{\infty}\frac{\left(  2-2^{-2n+s}\right)  }{\left(  2n+1\right)
!}\Gamma\left(  2n-s+1\right)   &  =\left[  -1+3^{s}+s\left(  2^{s}-2\right)
\right]  \Gamma\left(  -s\right) \nonumber\\
&  +\sqrt{\pi}\frac{\Gamma\left(  1-s\right)  \Gamma\left(  s\right)  }%
{\Gamma\left(  1+\frac{s}{2}\right)  \Gamma\left(  \frac{1+s}{2}\right)  }.
\end{align}
Using the above relation we can obtain the following functional relation%
\begin{align}
\zeta\left(  s\right)   &  =\frac{\pi^{s-1}}{\left(  1-2^{1-s}\right)  }%
\sin\left(  \frac{\pi s}{2}\right)  \left(  \sum_{n=1}^{\infty}\frac{\left(
2-2^{-2n+s}\right)  }{\left(  2n+1\right)  !}\Gamma\left(  2n-s+1\right)
\left[  \zeta\left(  2n-s+1\right)  -1\right]  \right. \nonumber\\
&  \left.  +\left[  -1+3^{s}+s\left(  2^{s}-2\right)  \right]  \Gamma\left(
-s\right)  +\sqrt{\pi}\frac{\Gamma\left(  1-s\right)  \Gamma\left(  s\right)
}{\Gamma\left(  1+\frac{s}{2}\right)  \Gamma\left(  \frac{1+s}{2}\right)
}\right)  . \label{newf2}%
\end{align}
Identity in Eq.~(\ref{newf2}) appears to be new.

\section{Conclusion}

In this work we have derived an integral representation for the Riemann zeta
function. Using this representation, we are able to obtain a new identity for
the Riemann zeta function. This representation allows us to express the
Riemann zeta function in terms of the Riemann zeta and Gamma functions. Note
that similar identities have been discovered by Milgram\cite{milgram}.

\end{document}